\documentclass{amsart}
\usepackage{amsfonts}

\newtheorem{theorem}{Theorem}
\theoremstyle{plain}

\newtheorem{corollary}{Corollary}

\newtheorem{lemma}{Lemma}

\numberwithin{equation}{section}
\input{tcilatex}

\begin{document}
\title[Inequalities for Functions of Bounded Variation]{Some Inequalities
for Functions of Bounded Variation with Applications to Landau Type Results}
\author{S.S. Dragomir}
\address{School of Computer Science and Mathematics\\
Victoria University of Technology\\
PO Box 14428\\
Melbourne VIC 8001, Australia.}
\email{sever.dragomir@vu.edu.au}
\urladdr{http://rgmia.vu.edu.au/SSDragomirWeb.html}
\date{06 August, 2004}

\begin{abstract}
Some inequalities for functions of bounded variation that provide reverses
for the inequality between the integral mean and the $p-$norm for $p\in %
\left[ 1,\infty \right] $ are established. Applications related to the
celebrated Landau inequality between the norms of the derivatives of a
function are also pointed out.
\end{abstract}

\maketitle

\section{Introduction}

The following inequality holding on finite intervals is well known.

If $f:\left[ a,b\right] \rightarrow \mathbb{R}$ is essentially bounded, then 
$f$ is integrable on $\left[ a,b\right] $ and 
\begin{equation}
\frac{1}{b-a}\left\vert \int_{a}^{b}f\left( t\right) dt\right\vert \leq
\left\Vert f\right\Vert _{\left[ a,b\right] ,\infty }  \label{1.1}
\end{equation}%
where $\left\Vert f\right\Vert _{\left[ a,b\right] ,\infty }:=ess\sup_{t\in %
\left[ a,b\right] }\left\vert f\left( t\right) \right\vert .$

The corresponding version in terms of $p-$norms, is the following H\"{o}lder
type inequality%
\begin{equation}
\frac{1}{\left( b-a\right) ^{1-\frac{1}{p}}}\left\vert \int_{a}^{b}f\left(
t\right) dt\right\vert \leq \left\Vert f\right\Vert _{\left[ a,b\right]
,p},\qquad p\geq 1,  \label{1.2}
\end{equation}%
provided $f\in L_{p}\left[ a,b\right] ,$ where%
\begin{equation*}
\left\Vert f\right\Vert _{\left[ a,b\right] ,p}:=\left(
\int_{a}^{b}\left\vert f\left( t\right) \right\vert ^{p}dt\right) ^{\frac{1}{%
p}},\qquad p\geq 1.
\end{equation*}

In the first part of this paper we point out some reverse inequalities for (%
\ref{1.1}) and (\ref{1.2}) in the case of functions of bounded variation.
These results are then employed in obtaining some Landau type inequalities.

For the latter, recall that if $I=\mathbb{R}_{+}$ or $I=\mathbb{R}$ and if $%
f:I\rightarrow \mathbb{R}$ is twice differentiable with $f,f^{\prime \prime
}\in L_{p}\left( I\right) ,$ $p\in \left[ 1,\infty \right] ,$ then $%
f^{\prime }\in L_{p}\left( I\right) .$ Moreover, there exists a constant $%
C_{p}\left( I\right) >0$ independent of the function $f,$ such that%
\begin{equation}
\left\Vert f^{\prime }\right\Vert _{p,I}\leq C_{p}\left( I\right) \left\Vert
f\right\Vert _{p,I}^{\frac{1}{2}}\left\Vert f^{\prime \prime }\right\Vert
_{p,I}^{\frac{1}{2}},  \label{1.3}
\end{equation}%
where $\left\Vert \cdot \right\Vert _{p,I}$ is the $p-$norm on the interval $%
I.$

The investigation of such inequalities was initiated by E. Landau \cite{LAN}
in 1914. He considered the case $p=\infty $ and proved that%
\begin{equation}
C_{\infty }\left( \mathbb{R}_{+}\right) =2\text{ and }C_{\infty }\left( 
\mathbb{R}\right) =\sqrt{2},  \label{1.4}
\end{equation}%
are the best constant for which (\ref{1.3}) holds.

For some classical and recent results related to Landau inequality, see \cite%
{DI},\cite{DP} and \cite{HL}-\cite{NI}.

\section{Some Reverse Inequalities on Bounded Intervals}

The following result for functions of bounded variation holds.

\begin{theorem}
\label{t1}Let $f:\left[ a,b\right] \rightarrow \mathbb{R}$ be a function of
bounded variation on $\left[ a,b\right] .$ Then%
\begin{equation}
\left\Vert f\right\Vert _{\left[ a,b\right] ,\infty }\leq \frac{1}{b-a}%
\left\vert \int_{a}^{b}f\left( t\right) dt\right\vert
+\bigvee\nolimits_{a}^{b}\left( f\right) .  \label{2.1}
\end{equation}%
The multiplicative constant $1$ in front of $\bigvee\nolimits_{a}^{b}\left(
f\right) $ cannot be replaced by a smaller quantity.
\end{theorem}

\begin{proof}
We apply the following Ostrowski type inequality obtained by the author in 
\cite{DR1} (see also \cite{DR2}):%
\begin{equation}
\left\vert f\left( x\right) -\frac{1}{b-a}\int_{a}^{b}f\left( t\right)
dt\right\vert \leq \left[ \frac{1}{2}+\frac{\left\vert x-\frac{a+b}{2}%
\right\vert }{b-a}\right] \bigvee\nolimits_{a}^{b}\left( f\right)
\label{2.2}
\end{equation}%
for any $x\in \left[ a,b\right] .$ The constant $\frac{1}{2}$ is best
possible in the sense that it cannot be replaced by a smaller quantity.

Taking the supremum in (\ref{2.2}) over $x\in \left[ a,b\right] $, we get%
\begin{align}
\left\Vert f-\frac{1}{b-a}\int_{a}^{b}f\left( t\right) dt\right\Vert _{\left[
a,b\right] ,\infty }& \leq \sup_{x\in \left[ a,b\right] }\left[ \frac{1}{2}+%
\frac{\left\vert x-\frac{a+b}{2}\right\vert }{b-a}\right] \bigvee%
\nolimits_{a}^{b}\left( f\right)  \label{2.3} \\
& =\bigvee\nolimits_{a}^{b}\left( f\right) .  \notag
\end{align}%
Now, by the triangle inequality applied for the sup-norm $\left\Vert \cdot
\right\Vert _{\infty },$ we get%
\begin{align*}
\left\Vert f\right\Vert _{\left[ a,b\right] ,\infty }& \leq \left\Vert f-%
\frac{1}{b-a}\int_{a}^{b}f\left( t\right) dt\right\Vert _{\left[ a,b\right]
,\infty }+\left\vert \frac{1}{b-a}\int_{a}^{b}f\left( t\right) dt\right\vert
\\
& \leq \frac{1}{b-a}\left\vert \int_{a}^{b}f\left( t\right) dt\right\vert
+\bigvee\nolimits_{a}^{b}\left( f\right)
\end{align*}%
and the inequality (\ref{2.1}) is proved.

To prove the sharpness of the constant $1,$ assume that the following
inequality holds%
\begin{equation}
\left\Vert f\right\Vert _{\left[ a,b\right] ,\infty }\leq \frac{1}{b-a}%
\left\vert \int_{a}^{b}f\left( t\right) dt\right\vert
+C\bigvee\nolimits_{a}^{b}\left( f\right)  \label{2.4}
\end{equation}%
with a $C>0.$

Consider the function $f_{0}:\left[ a,b\right] \rightarrow \mathbb{R}$,%
\begin{equation*}
f_{0}\left( t\right) =\left\{ 
\begin{array}{ll}
0, & t\in \lbrack a,b) \\ 
&  \\ 
1, & t=b.%
\end{array}%
\right.
\end{equation*}%
Then $f_{0}$ is of bounded variation on $\left[ a,b\right] $ and%
\begin{equation*}
\left\Vert f_{0}\right\Vert _{\left[ a,b\right] ,\infty }=1,\ \ \ \
\int_{a}^{b}f_{0}\left( t\right) dt=0\text{ \ \ and \ \ }\bigvee%
\nolimits_{a}^{b}\left( f_{0}\right) =1.
\end{equation*}%
For this choice, (\ref{2.4}) becomes $C\geq 1,$ proving the sharpness of the
constant.
\end{proof}

The corresponding result for $p-$norms, where $p\geq 1,$ is embodied in the
following theorem.

\begin{theorem}
\label{t2}Let $f:\left[ a,b\right] \rightarrow \mathbb{R}$ be a function of
bounded variation on $\left[ a,b\right] .$ Then for $p\geq 1$ one has the
inequality%
\begin{equation}
\left\Vert f\right\Vert _{\left[ a,b\right] ,p}\leq \frac{1}{\left(
b-a\right) ^{1-\frac{1}{p}}}\left\vert \int_{a}^{b}f\left( t\right)
dt\right\vert +\frac{1}{2}\cdot \frac{\left( b-a\right) ^{\frac{1}{p}}\left(
2^{p+1}-1\right) ^{\frac{1}{p}}}{\left( p+1\right) ^{\frac{1}{p}}}%
\bigvee\nolimits_{a}^{b}\left( f\right) .  \label{2.5}
\end{equation}%
The constant $\frac{1}{2}$ is best possible in the sense that it cannot be
replaced by a smaller quantity.
\end{theorem}

\begin{proof}
Taking the $p-$norm in (\ref{2.2}), we deduce%
\begin{equation*}
\left\Vert f-\frac{1}{b-a}\int_{a}^{b}f\left( t\right) dt\right\Vert _{\left[
a,b\right] ,p}\leq \bigvee\nolimits_{a}^{b}\left( f\right) I_{p},
\end{equation*}%
where%
\begin{equation*}
I_{p}:=\left( \int_{a}^{b}\left[ \frac{1}{2}+\frac{\left\vert x-\frac{a+b}{2}%
\right\vert }{b-a}\right] ^{p}dx\right) ^{\frac{1}{p}},\ \ \ p\geq 1.
\end{equation*}%
We observe that%
\begin{align*}
I_{p}& :=\left( \int_{a}^{\frac{a+b}{2}}\left[ \frac{1}{2}+\frac{\frac{a+b}{2%
}-x}{b-a}\right] ^{p}dx+\int_{\frac{a+b}{2}}^{b}\left[ \frac{1}{2}+\frac{x-%
\frac{a+b}{2}}{b-a}\right] ^{p}dx\right) ^{\frac{1}{p}} \\
& =\frac{1}{b-a}\left[ \int_{a}^{\frac{a+b}{2}}\left( b-x\right)
^{p}dx+\int_{\frac{a+b}{2}}^{b}\left( x-a\right) ^{p}dx\right] \\
& =\frac{\left( b-a\right) ^{\frac{1}{p}}\left( 2^{p+1}-1\right) ^{\frac{1}{p%
}}}{2\left( p+1\right) ^{\frac{1}{p}}},\ \ \ p\geq 1.
\end{align*}%
Using the triangle inequality for the $p-$norm $\left\Vert \cdot \right\Vert
_{p},$ we get%
\begin{align*}
\left\Vert f\right\Vert _{\left[ a,b\right] ,p}& \leq \left\Vert f-\frac{1}{%
b-a}\int_{a}^{b}f\left( t\right) dt\right\Vert _{\left[ a,b\right]
,p}+\left\Vert \frac{1}{b-a}\int_{a}^{b}f\left( t\right) dt\right\Vert _{%
\left[ a,b\right] ,p} \\
& \leq \frac{\left( b-a\right) ^{\frac{1}{p}}\left( 2^{p+1}-1\right) ^{\frac{%
1}{p}}}{2\left( p+1\right) ^{\frac{1}{p}}}\bigvee\nolimits_{a}^{b}\left(
f\right) +\left( b-a\right) ^{\frac{1}{p}}\left\vert \frac{1}{b-a}%
\int_{a}^{b}f\left( t\right) dt\right\vert
\end{align*}%
and the inequality (\ref{2.5}) is obtained.

Now, assume that (\ref{2.5}) holds with a constant ~$D>0$ instead of $\frac{1%
}{2},$ i.e.,%
\begin{equation}
\left\Vert f\right\Vert _{\left[ a,b\right] ,p}\leq \frac{1}{\left(
b-a\right) ^{1-\frac{1}{p}}}\left\vert \int_{a}^{b}f\left( t\right)
dt\right\vert +D\cdot \frac{\left( b-a\right) ^{\frac{1}{p}}\left(
2^{p+1}-1\right) ^{\frac{1}{p}}}{\left( p+1\right) ^{\frac{1}{p}}}%
\bigvee\nolimits_{a}^{b}\left( f\right) .  \label{2.6}
\end{equation}%
Consider the function $f_{0}:\left[ a,b\right] \rightarrow \mathbb{R}$ with $%
a=0$ and $b>1$ given by%
\begin{equation*}
f_{0}\left( t\right) =\left\{ 
\begin{array}{ll}
0, & \text{if \ }t\in \left[ 0,b-1\right]  \\ 
&  \\ 
1, & \text{if \ }t\in (b-1,b].%
\end{array}%
\right. 
\end{equation*}%
This function is of bounded variation on $\left[ a,b\right] $ and%
\begin{equation*}
\left\Vert f\right\Vert _{\left[ a,b\right] ,p}=1,\ \ \ \
\int_{a}^{b}f\left( t\right) dt=1\text{ \ \ and \ \ }\bigvee%
\nolimits_{a}^{b}\left( f\right) =1
\end{equation*}%
and then, by (\ref{2.6}), we deduce%
\begin{equation*}
1\leq \frac{1}{b^{1-\frac{1}{p}}}+D\frac{b^{\frac{1}{p}}\left(
2^{p+1}-1\right) ^{\frac{1}{p}}}{\left( p+1\right) ^{\frac{1}{p}}},\ \ \
b>1,\ p\geq 1
\end{equation*}%
giving%
\begin{equation}
b^{1-\frac{1}{p}}\leq 1+D\cdot b\frac{\left( 2^{p+1}-1\right) ^{\frac{1}{p}}%
}{\left( p+1\right) ^{\frac{1}{p}}}.  \label{2.7}
\end{equation}%
Denote%
\begin{equation*}
q:=\frac{\left( 2^{p+1}-1\right) ^{\frac{1}{p}}}{\left( p+1\right) ^{\frac{1%
}{p}}}.
\end{equation*}%
Then%
\begin{equation*}
\ln q=\frac{\ln \left( 2^{p+1}-1\right) -\ln \left( p+1\right) }{p}.
\end{equation*}%
We observe, by L'Hospital theorem that%
\begin{align*}
\lim_{p\rightarrow \infty }\left[ \frac{\ln \left( 2^{p+1}-1\right) }{p}%
\right] & =\lim_{p\rightarrow \infty }\frac{\left[ \ln \left(
2^{p+1}-1\right) \right] ^{\prime }}{\left( p\right) ^{\prime }} \\
& =\lim_{p\rightarrow \infty }\frac{\left( 2^{p+1}-1\right) ^{\prime }}{%
2^{p+1}-1}=\ln 2
\end{align*}%
and%
\begin{equation*}
\lim_{p\rightarrow \infty }\left[ \frac{\ln \left( p+1\right) }{p}\right] =0,
\end{equation*}%
consequently%
\begin{equation*}
\lim_{p\rightarrow \infty }q=2.
\end{equation*}%
Taking the limit over $p\rightarrow \infty $ in (\ref{2.7}), we deduce%
\begin{equation*}
b\leq 1+2Db,\text{ \ \ for \ }b>1
\end{equation*}%
from where we get%
\begin{equation}
D\geq \frac{b-1}{2b},\ \ \ b>1.  \label{2.8}
\end{equation}%
Taking the limit over $b\rightarrow \infty $ in (\ref{2.8}) we conclude that 
$D\geq \frac{1}{2},$ showing that the constant $\frac{1}{2}$ in (\ref{2.5})
cannot be replaced by a smaller quantity in (\ref{2.5}).
\end{proof}

\section{Some Inequalities of Landau Type on Unbounded Intervals}

The following technical lemma will be used in the following (see also \cite%
{DP}).

\begin{lemma}
\label{l3.1}Let $C,D>0$ and $r,u\in (0,1].$ Consider the function $%
g_{r,u}:\left( 0,\infty \right) \rightarrow \left( 0,\infty \right) $ given
by%
\begin{equation}
g_{r,u}\left( \lambda \right) =\frac{C}{\lambda ^{u}}+D\lambda ^{r}.
\label{3.1}
\end{equation}%
Define%
\begin{equation}
\lambda _{0}:=\left( \frac{uC}{rD}\right) ^{\frac{1}{r+u}}\in \left(
0,\infty \right) .  \label{3.2}
\end{equation}%
Then we have%
\begin{equation}
\inf_{\lambda \in \left( 0,\infty \right) }g_{r,u}\left( \lambda \right)
=g\left( \lambda _{0}\right) =\frac{r+u}{u^{\frac{u}{r+u}}\cdot r^{\frac{u}{%
r+u}}}C^{\frac{r}{r+u}}D^{\frac{r}{r+u}}.  \label{3.3}
\end{equation}
\end{lemma}

\begin{proof}
We observe that%
\begin{equation*}
g_{r,u}^{\prime }\left( \lambda \right) =\frac{rD\lambda ^{r+u}-Cu}{\lambda
^{u+1}},\ \ \ \lambda \in \left( 0,\infty \right) .
\end{equation*}%
The unique solution of the equation $g_{r,u}^{\prime }\left( \lambda \right)
=0,$ $\lambda \in \left( 0,\infty \right) $ is $\lambda _{0}$ provided by (%
\ref{3.2}).

The function $g_{r,u}$ is decreasing on $\left( 0,\lambda _{0}\right) $ and
increasing on $\left( \lambda _{0},\infty \right) .$ The global minimum for $%
g_{r,u}$ on $\left( 0,\infty \right) $ is%
\begin{align*}
g_{r,u}\left( \lambda _{0}\right) & =\frac{C}{\left( \frac{uC}{rD}\right) ^{%
\frac{u}{r+u}}}+D\left( \frac{uC}{rD}\right) ^{\frac{r}{r+u}} \\
& =\frac{r+u}{u^{\frac{u}{r+u}}r^{\frac{r}{r+u}}}C^{\frac{r}{r+u}}D^{\frac{u%
}{r+u}}
\end{align*}%
and the equality (\ref{3.3}) is proved.
\end{proof}

The following particular cases are useful in applications.

\begin{corollary}
\label{c3.2}Let $C,D>0.$

\begin{enumerate}
\item[(i)] For $r\in (0,1],$ consider the function $g_{r}:\left( 0,\infty
\right) \rightarrow \left( 0,\infty \right) ,$ given by%
\begin{equation}
g_{r}\left( \lambda \right) =\frac{C}{\lambda }+D\lambda ^{r}.  \label{3.4}
\end{equation}%
Define%
\begin{equation}
\overline{\lambda _{0}}=\left( \frac{C}{rD}\right) ^{\frac{1}{r+1}}\in
\left( 0,\infty \right) .  \label{3.5}
\end{equation}%
Then we have%
\begin{equation}
\inf_{\lambda \in \left( 0,\infty \right) }g_{r}\left( \lambda \right)
=g_{r}\left( \overline{\lambda _{0}}\right) =\frac{r+1}{r^{\frac{r}{r+u}}}C^{%
\frac{r}{r+1}}D^{\frac{1}{r+1}}.  \label{3.6}
\end{equation}

\item[(ii)] For $u\in (0,1],$ consider the function $g_{u}:\left( 0,\infty
\right) \rightarrow \left( 0,\infty \right) ,$ given by%
\begin{equation}
g_{u}\left( \lambda \right) =\frac{C}{\lambda ^{u}}+D\lambda .  \label{3.7}
\end{equation}%
Define%
\begin{equation*}
\widetilde{\lambda _{0}}=\left( \frac{uC}{D}\right) ^{\frac{1}{1+u}}\in
\left( 0,\infty \right) .
\end{equation*}%
Then we have%
\begin{equation*}
\inf_{\lambda \in \left( 0,\infty \right) }g_{u}\left( \lambda \right)
=g_{u}\left( \widetilde{\lambda _{0}}\right) =\frac{1+u}{u^{\frac{u}{1+u}}}%
C^{\frac{1}{u+1}}D^{\frac{u}{u+1}}.
\end{equation*}
\end{enumerate}
\end{corollary}

The following result holds.

\begin{theorem}
\label{t3.3}Let $J$ be an unbounded subinterval of $\mathbb{R}$ and $%
g:J\rightarrow \mathbb{R}$ a locally absolutely continuous function on $J.$
If $g\in L_{\infty }\left( J\right) $, the derivative $g^{\prime
}:J\rightarrow \mathbb{R}$ is of locally bounded variation and there exists
a constant $V_{J}>0$ and $r\in (0,1]$ such that%
\begin{equation}
\left\vert \bigvee\nolimits_{a}^{b}\left( g^{\prime }\right) \right\vert
\leq V_{J}\left\vert a-b\right\vert ^{r}\text{ \ \ for any \ }a,b\in J;
\label{3.8}
\end{equation}%
then $g^{\prime }\in L_{\infty }\left( J\right) $ and one has the inequality%
\begin{equation}
\left\Vert g^{\prime }\right\Vert _{J,\infty }\leq \frac{2^{\frac{r}{r+1}%
}\left( r+1\right) }{r^{\frac{r}{r+1}}}\left\Vert g\right\Vert _{J,\infty }^{%
\frac{r}{r+1}}V_{J}^{\frac{1}{r+1}}.  \label{3.9}
\end{equation}
\end{theorem}

\begin{proof}
Applying Theorem \ref{t1} for the function $f=g^{\prime }$ on $\left[ a,b%
\right] $ (or $\left[ b,a\right] $), we deduce%
\begin{equation}
\left\Vert g^{\prime }\right\Vert _{\left[ a,b\right] ,\infty }\leq \frac{%
\left\vert g\left( b\right) -g\left( a\right) \right\vert }{\left\vert
b-a\right\vert }+\left\vert \bigvee\nolimits_{a}^{b}\left( g^{\prime
}\right) \right\vert .  \label{3.10}
\end{equation}%
for any $a,b\in J,$ $a\neq b.$

Since $\left\vert g^{\prime }\left( b\right) \right\vert \leq \left\Vert
g^{\prime }\right\Vert _{\left[ a,b\right] ,\infty },$ $\left\vert g\left(
b\right) -g\left( a\right) \right\vert \leq 2\left\Vert g\right\Vert
_{J,\infty },$ then by (\ref{3.8}) and (\ref{3.10}) we deduce%
\begin{equation}
\left\vert g^{\prime }\left( b\right) \right\vert \leq \frac{2\left\Vert
g\right\Vert _{J,\infty }}{\left\vert b-a\right\vert }+V_{J}\left\vert
b-a\right\vert ^{r}  \label{3.11}
\end{equation}%
for any $a,b\in J,$ $a\neq b.$

Fix $b\in J.$ Then for any $\lambda >0,$ there exists an $a\in J$ such that $%
\lambda =\left\vert b-a\right\vert .$ Consequently, by (\ref{3.11}), we
deduce that%
\begin{equation}
\left\vert g^{\prime }\left( b\right) \right\vert \leq \frac{2\left\Vert
g\right\Vert _{J,\infty }}{\lambda }+V_{J}\lambda ^{r}  \label{3.12}
\end{equation}%
for any $\lambda >0$ and $b\in J.$

Taking the infimum over $\lambda \in \left( 0,\infty \right) $ in (\ref{3.12}%
) and using Corollary \ref{c3.2}, we deduce%
\begin{align}
\left\vert g^{\prime }\left( b\right) \right\vert & \leq \frac{r+1}{r^{\frac{%
r}{r+1}}}\left( 2\left\Vert g\right\Vert _{J,\infty }\right) ^{\frac{r}{r+1}%
}\cdot V_{J}^{\frac{1}{r+1}}  \label{3.13} \\
& =\frac{2^{\frac{r}{r+1}}\left( r+1\right) }{r^{\frac{r}{r+1}}}\left\Vert
g\right\Vert _{J,\infty }^{\frac{r}{r+1}}V_{J}^{\frac{1}{r+1}}  \notag
\end{align}%
for any $b\in J.$ Finally, taking the supremum in (\ref{3.13}) over $b\in J,$
we deduce the desired result (\ref{3.9}).
\end{proof}

There are a number of particular cases of interest.

\begin{corollary}
\label{c3.4}Assume that $g:J\rightarrow \mathbb{R}$ is such that $g^{\prime
}:J\rightarrow \mathbb{R}$ is locally absolutely continuous and $g^{\prime
\prime }\in L_{\infty }\left( J\right) .$ If $g\in L_{\infty }\left(
J\right) ,$ then $g^{\prime }\in L_{\infty }\left( J\right) $ and%
\begin{equation}
\left\Vert g^{\prime }\right\Vert _{J,\infty }\leq 2\sqrt{2}\left\Vert
g\right\Vert _{J,\infty }^{\frac{1}{2}}\left\Vert g^{\prime \prime
}\right\Vert _{J,\infty }^{\frac{1}{2}}.  \label{3.14}
\end{equation}
\end{corollary}

\begin{proof}
If $g^{\prime \prime }\in L_{\infty }\left( J\right) ,$ then%
\begin{equation*}
\left\vert \bigvee\nolimits_{a}^{b}\left( g^{\prime }\right) \right\vert
=\left\vert \int_{a}^{b}\left\vert g^{\prime \prime }\left( t\right)
\right\vert dt\right\vert \leq \left\vert b-a\right\vert \left\Vert
g^{\prime \prime }\right\Vert _{J,\infty }
\end{equation*}%
for any $a,b\in J,$ giving, by (\ref{3.11}), that%
\begin{equation}
\left\vert g^{\prime }\left( b\right) \right\vert \leq \frac{2\left\Vert
g\right\Vert _{J,\infty }}{\left\vert b-a\right\vert }+\left\Vert g^{\prime
\prime }\right\Vert _{J,\infty }\left\vert b-a\right\vert  \label{3.15}
\end{equation}%
for any $a,b\in J,$ $a\neq b.$

Applying Theorem \ref{t3.3} for $V_{J}=\left\Vert g^{\prime \prime
}\right\Vert _{J,\infty }$ and $r=1,$ we deduce (\ref{3.14}).
\end{proof}

The following result is also of interest.

\begin{corollary}
\label{c3.5}Assume that $g:J\rightarrow \mathbb{R}$ is such that $g^{\prime
}\in L_{p}\left( J\right) ,$ $p>1.$ If $g\in L_{\infty }\left( J\right) ,$
then $g^{\prime }\in L_{\infty }\left( J\right) $ and%
\begin{equation}
\left\Vert g^{\prime }\right\Vert _{J,\infty }\leq \frac{2^{\frac{p-1}{2p-1}%
}\left( 2p-1\right) }{\left( p-1\right) ^{\frac{p-1}{2p-1}}p^{\frac{p}{2p-1}}%
}\cdot \left\Vert g\right\Vert _{J,\infty }^{\frac{p-1}{2p-1}}\left\Vert
g^{\prime \prime }\right\Vert _{J,p}^{\frac{p-1}{2p-1}}.  \label{3.16}
\end{equation}
\end{corollary}

\begin{proof}
Using H\"{o}lder's inequality, we have%
\begin{align*}
\left\vert \bigvee\nolimits_{a}^{b}\left( g^{\prime }\right) \right\vert &
=\left\vert \int_{a}^{b}\left\vert g^{\prime \prime }\left( t\right)
\right\vert dt\right\vert \leq \left\vert \int_{a}^{b}dt\right\vert ^{\frac{1%
}{q}}\left\vert \int_{a}^{b}\left\vert g^{\prime \prime }\left( t\right)
\right\vert ^{p}dt\right\vert ^{\frac{1}{p}} \\
& \leq \left\vert b-a\right\vert ^{\frac{1}{q}}\left\Vert g^{\prime \prime
}\right\Vert _{J,p},\ \ \ p>1,\ \frac{1}{p}+\frac{1}{q}=1
\end{align*}%
for any $a,b\in J$, giving, by (\ref{3.11}), that%
\begin{equation}
\left\vert g^{\prime }\left( b\right) \right\vert \leq \frac{2\left\Vert
g\right\Vert _{J,\infty }}{\left\vert b-a\right\vert }+\left\vert
b-a\right\vert ^{\frac{1}{q}}\left\Vert g^{\prime \prime }\right\Vert _{J,p},
\label{3.17}
\end{equation}%
for any $a,b\in J,$ $a\neq b.$

Applying Theorem \ref{t3.3} for $V_{J}=\left\Vert g^{\prime \prime
}\right\Vert _{J,p}$ and $r=\frac{1}{q}=\frac{p-1}{p},$ we deduce (\ref{3.16}%
).
\end{proof}

The following result also holds.

\begin{theorem}
\label{t3.6}Let $J$ be an unbounded subinterval of $\mathbb{R}$ and $g:J%
\mathbb{\rightarrow R}$ a locally absolutely continuous function on $J.$ If $%
g^{\prime }\in L_{1}\left( J\right) ,$ the derivative $g^{\prime
}:J\rightarrow \mathbb{R}$ is of locally bounded variation and there exists
a constant $V_{J}>0$ and $r\in (0,1]$ such that%
\begin{equation}
\left\vert \bigvee\nolimits_{a}^{b}\left( g^{\prime }\right) \right\vert
\leq V_{J}\left\vert a-b\right\vert ^{r}\text{ \ \ for any \ }a,b\in J;
\label{3.18}
\end{equation}%
then $g^{\prime }\in L_{\infty }\left( J\right) $ and one has the inequality%
\begin{equation}
\left\Vert g^{\prime }\right\Vert _{J,\infty }\leq \frac{r+1}{r^{\frac{r}{r+1%
}}}\left\Vert g^{\prime }\right\Vert _{J,1}^{\frac{r}{r+1}}V_{J}^{\frac{1}{%
r+1}}.  \label{3.19}
\end{equation}
\end{theorem}

\begin{proof}
Since, for any $a,b\in J,$%
\begin{equation*}
\left\vert g\left( b\right) -g\left( a\right) \right\vert \leq \left\vert
\int_{a}^{b}g^{\prime }\left( s\right) ds\right\vert \leq \left\vert
\int_{a}^{b}\left\vert g^{\prime }\left( s\right) \right\vert ds\right\vert
\leq \left\Vert g^{\prime }\right\Vert _{J,1},
\end{equation*}%
then, by (\ref{3.10}) and (\ref{3.18}), we deduce%
\begin{equation*}
\left\vert g^{\prime }\left( b\right) \right\vert \leq \frac{\left\Vert
g^{\prime }\right\Vert _{J,1}}{\left\vert b-a\right\vert }+V_{J}\left\vert
b-a\right\vert ^{r}
\end{equation*}%
for any $a,b\in J,$ $a\neq b.$

Using an argument similar to the one in Theorem \ref{t3.3}, we deduce (\ref%
{3.19}).
\end{proof}

The following particular case also holds.

\begin{corollary}
\label{c3.7}Assume that $g:J\rightarrow \mathbb{R}$ is such that $g^{\prime
}:J\rightarrow \mathbb{R}$ is locally absolutely continuous and $g^{\prime
\prime }\in L_{\infty }\left( J\right) .$ If $g^{\prime }\in L_{1}\left(
J\right) ,$ then $g^{\prime }\in L_{\infty }\left( J\right) $ and%
\begin{equation}
\left\Vert g^{\prime }\right\Vert _{J,\infty }\leq 2\left\Vert g^{\prime
}\right\Vert _{J,1}^{\frac{1}{2}}\left\Vert g^{\prime \prime }\right\Vert
_{J,\infty }^{\frac{1}{2}}.  \label{3.20}
\end{equation}
\end{corollary}

\begin{corollary}
\label{c3.8}Assume that $g:J\rightarrow \mathbb{R}$ is such that $g^{\prime
}:J\rightarrow \mathbb{R}$ is locally absolutely continuous and $g^{\prime
\prime }\in L_{p}\left( J\right) ,$ $p>1.$ If $g^{\prime }\in L_{1}\left(
J\right) ,$ then $g^{\prime }\in L_{\infty }\left( J\right) $ and%
\begin{equation}
\left\Vert g^{\prime }\right\Vert _{J,\infty }\leq \frac{2p-1}{\left(
p-1\right) ^{\frac{p-1}{2p-1}}p^{\frac{p}{2p-1}}}\cdot \left\Vert
g\right\Vert _{J,1}^{\frac{p-1}{2p-1}}\left\Vert g^{\prime \prime
}\right\Vert _{J,p}^{\frac{p-1}{2p-1}}.  \label{3.21}
\end{equation}
\end{corollary}

We may state the following result as well.

\begin{theorem}
\label{t3.9}Let $J$ be an unbounded subinterval of $\mathbb{R}$ and $%
g:J\rightarrow \mathbb{R}$ a locally absolutely continuous function on $J.$
If $g^{\prime }\in L_{\alpha }\left( J\right) ,$ $\alpha >1,$ the derivative 
$g^{\prime }:J\rightarrow \mathbb{R}$ is of locally bounded variation on $J$
and there exists a constant $V_{J}>0$ and $r\in (0,1]$ such that%
\begin{equation}
\left\vert \bigvee\nolimits_{a}^{b}\left( g^{\prime }\right) \right\vert
\leq V_{J}\left\vert b-a\right\vert ^{r}\text{ \ \ for any \ }a,b\in J;
\label{3.22}
\end{equation}%
then $g^{\prime }\in L_{\infty }\left( J\right) $ and one has the inequality%
\begin{equation}
\left\Vert g^{\prime }\right\Vert _{J,\infty }\leq \frac{\alpha r+1}{\alpha
^{\frac{\alpha r}{\alpha r+1}}r^{\frac{\alpha r}{\alpha r+1}}}\left\Vert
g^{\prime }\right\Vert _{J,\alpha }^{\frac{\alpha r}{\alpha r+1}}V_{J}^{%
\frac{1}{\alpha r+1}}.  \label{3.23}
\end{equation}
\end{theorem}

\begin{proof}
By H\"{o}lder's integral inequality, we have%
\begin{align*}
\left\vert g\left( b\right) -g\left( a\right) \right\vert & =\left\vert
\int_{a}^{b}g^{\prime }\left( s\right) ds\right\vert \leq \left\vert
\int_{a}^{b}\left\vert g^{\prime }\left( s\right) \right\vert ds\right\vert
\\
& \leq \left\vert b-a\right\vert ^{\frac{1}{\beta }}\left\Vert g^{\prime
}\right\Vert _{J,\alpha },\ \ \ \alpha >1,\ \frac{1}{\alpha }+\frac{1}{\beta 
}=1,
\end{align*}%
and then, by (\ref{3.10}) and (\ref{3.18}), we deduce%
\begin{align}
\left\vert g^{\prime }\left( b\right) \right\vert & \leq \frac{\left\vert
b-a\right\vert ^{\frac{1}{\beta }}\left\Vert g^{\prime }\right\Vert
_{J,\alpha }}{\left\vert b-a\right\vert }+\left\vert b-a\right\vert ^{r}V_{J}
\label{3.24} \\
& =\frac{\left\Vert g^{\prime }\right\Vert _{J,\alpha }}{\left\vert
b-a\right\vert ^{\frac{1}{\alpha }}}+\left\vert b-a\right\vert ^{r}V_{J} 
\notag
\end{align}%
for any $a,b\in J,$ $a\neq b.$

Fix $b\in J.$ Then for any $\lambda >0,$ there exists an $a\in J$ such that $%
\lambda =\left\vert b-a\right\vert .$ Consequently, by (\ref{3.14}) we
deduce that%
\begin{equation}
\left\vert g^{\prime }\left( b\right) \right\vert \leq \frac{\left\Vert
g^{\prime }\right\Vert _{J,\alpha }}{\lambda ^{\frac{1}{\alpha }}}+\lambda
^{r}V_{J}  \label{3.25}
\end{equation}%
for any $\lambda >0$ and $b\in J.$

Taking the infimum over $\lambda \in \left( 0,\infty \right) $ in (\ref{3.25}%
) and using Lemma \ref{l3.1} for $u=\frac{1}{\alpha },$ we deduce%
\begin{align*}
\left\vert g^{\prime }\left( b\right) \right\vert & \leq \frac{r+\frac{1}{%
\alpha }}{\left( \frac{1}{\alpha }\right) ^{\frac{\frac{1}{\alpha }}{r+\frac{%
1}{\alpha }}}r^{\frac{r}{r+\frac{1}{\alpha }}}}\left\Vert g^{\prime
}\right\Vert _{J,\alpha }^{\frac{r}{r+\frac{1}{\alpha }}}V_{J}^{\frac{\frac{1%
}{\alpha }}{r+\frac{1}{\alpha }}} \\
& =\frac{\alpha r+1}{\alpha ^{\frac{\alpha r}{\alpha r+1}}r^{\frac{\alpha r}{%
\alpha r+1}}}\left\Vert g^{\prime }\right\Vert _{J,\alpha }^{\frac{\alpha r}{%
\alpha r+1}}V_{J}^{\frac{1}{\alpha r+1}}
\end{align*}%
for any $b\in J,$ giving the desired result (\ref{3.23}).
\end{proof}

The following corollary holds.

\begin{corollary}
\label{c3.9}Assume that $g:J\rightarrow \mathbb{R}$ is such that $g^{\prime
} $ is locally absolutely continuous and $g^{\prime \prime }\in L_{\infty
}\left( J\right) .$ If $g^{\prime }\in L_{\alpha }\left( J\right) ,$ $\alpha
>1,$ then $g^{\prime }\in L_{\infty }\left( J\right) $ and%
\begin{equation}
\left\Vert g^{\prime }\right\Vert _{J,\infty }\leq \frac{\alpha +1}{\alpha ^{%
\frac{\alpha }{\alpha +1}}r}\left\Vert g^{\prime }\right\Vert _{J,\alpha }^{%
\frac{\alpha }{\alpha +1}}\left\Vert g^{\prime \prime }\right\Vert
_{J,\infty }^{\frac{1}{\alpha +1}}.  \label{3.26}
\end{equation}
\end{corollary}

Finally we have

\begin{corollary}
\label{c3.10}Assume that $g:J\rightarrow \mathbb{R}$ is such that $g^{\prime
}$ is locally absolutely continuous and $g^{\prime \prime }\in L_{p}\left(
J\right) ,$ $p>1.$ If $g^{\prime }\in L_{\alpha }\left( J\right) ,$ $\alpha
>1,$ then $g^{\prime }\in L_{\infty }\left( J\right) $ and%
\begin{equation}
\left\Vert g^{\prime }\right\Vert _{J,\infty }\leq \frac{\alpha \left(
p-1\right) +p}{\alpha ^{\frac{\alpha \left( p-1\right) }{\alpha \left(
p-1\right) +p}}\left( p-1\right) ^{\frac{\alpha \left( p-1\right) }{\alpha
\left( p-1\right) +p}}\cdot p^{\frac{p}{\alpha \left( p-1\right) +p}}}%
\left\Vert g^{\prime }\right\Vert _{J,\alpha }^{\frac{\alpha \left(
p-1\right) }{\alpha \left( p-1\right) +p}}\left\Vert g^{\prime \prime
}\right\Vert _{J,p}^{\frac{p}{\alpha \left( p-1\right) +p}}.  \label{3.27}
\end{equation}
\end{corollary}

\end{document}